\documentclass[12pt,a4paper]{article}
\usepackage{amsmath}
\usepackage{amsfonts}
\usepackage{amssymb}
\usepackage{mathrsfs}
\usepackage{amsthm}
\usepackage{mathrsfs}
\usepackage{bm}

\usepackage{framed}
\usepackage{color}

\usepackage{ascmac}

\usepackage{latexsym}

\usepackage[top=15truemm,bottom=25truemm,left=25truemm,right=25truemm]{geometry}

\definecolor{lightgray}{rgb}{0.75,0.75,0.75}

\title{Random discretization of O'Hara knot energy}
\author{Jun Okamoto\thanks{Graduate School of Mathematical Sciences, The University of Tokyo, Komaba 3-8-1 Meguro Tokyo 153-8914, Japan. E-mail: \texttt{okamoto@ms.u-tokyo.ac.jp}}}
\date{May 16, 2019}

\newtheorem{Thm}{Theorem}[section]
\newtheorem{Prop}{Proposition}[section]

\theoremstyle{definition}
\newtheorem{Dfn}{Definition}[section]
\newtheorem{Rem}{Remark}

\newcommand{\R}{\mathbb{R}}

\newcommand{\N}{\mathbb{N}}

\newcommand{\T}{\mathbb{R}/\mathbb{Z}}
\newcommand{\TL}{\mathbb{R}/\mathcal{L}\mathbb{Z}}

\newcommand{\D}{\displaystyle}

\newcommand{\ep}{\varepsilon}
\newcommand{\al}{\alpha}

\newcommand{\relmiddle}[1]{\mathrel{}\middle#1\mathrel{}}

\newcounter{mynum}

\makeatletter

\begin{document}
\maketitle

\begin{abstract}

In this paper, we considered random discrete approximation of
O'Hara energy. O'Hara energy is
the energy defined for a knot, and O'Hara energy
was introduced for defining the standard
shape for each knot class (equivalence class by
ambient isotopy) by variational method. In the
case of a specific exponent, due to energy invariance
under M\"{o}bius transformation, this energy
is called M\"{o}bius energy. Although discretization
for various M\"{o}bius energies has been defined to analyse the shape of the minimizer so
far, only $\Gamma$-convergence to the original energy
has been shown for a conventional discretization.
In this study, we are successful to show
locally uniform convergence and compactness
of discrete energy in a space based on optimal
transport theory, by introducing random
discrete approximation of O'Hara energy using
random variable.
\end{abstract}
\section{Introduction}

Let $\mathcal{A}$ be a set of all closed regular curve that is parametrised by arc length in $\R^n,$ with no self-intersections, and with total length $\mathcal{L}$ i.e. $\mathcal{A}:=\{\gamma \in C^{0,1}(\TL,\R^n) \mid |\gamma'(x)|=1 \ \ \mbox{a.e.}~x \in \TL \},$ and $\al,p \in (0,\infty),$ the O'Hara $(\alpha,p)$-energy $\mathcal{E}^{\alpha,p}:\mathcal{A}\rightarrow \R \cup \{+\infty\}$ is defined as follow:
\begin{align}
\D \mathcal{E}^{\al,p}(\gamma) :=\int_{(\mathbb{R}/\mathcal{L}\mathbb{Z})^2}  \mathcal{M}^{\al,p}(\gamma) dx dy,
\end{align}
where
\begin{align}
\mathcal{M}^{\al,p}(\gamma) (x,y)=\left(\frac{1}{|\gamma(x)-\gamma(y) |^\al}-\frac{1}{\mathcal{D}(\gamma(x),\gamma(y))^\al} \right)^p  \qquad (x,y) \in (\TL)^2
\end{align} and $\mathcal{D}$ is the length of shortest arc of the curve $\gamma$ connecting  the two points $\gamma(x)$ and $\gamma(y),$ i.e.
\begin{align}
\mathcal{D}(\gamma(x),\gamma(y))=\min \{\mathcal{L}-|x-y|,|x-y|\}.
\end{align}

This energy introduced and investigated by O'Hara in \cite{O1}-\cite{O4} for defining the standard shape for each knot class by variational method.

In the case of $\al=2,p=1$, due to energy invariance under M\"{o}bius transformation, this energy called ''M\"{o}bius energy". It is possible to show the existence of minimizer in the ''prime knot". R. Kusner and J. Sullivan conjectured the minimizer in composite knot class may not exist \cite{Ku}. This conjecture was established by numerical calculation with discretization of M\"{o}bius energy.

In this paper, we introduce the weighted O'Hara energy $\mathcal{E}^{\al,p}_\rho:\mathcal{A}\rightarrow \R \cup \{+\infty\}$ with weight $\rho:\TL \rightarrow \R_{\ge 0}$ defined
\begin{align}
\mathcal{E}^{\al,p}_\rho(\gamma) := \int_{(\mathbb{R}/\mathcal{L}\mathbb{Z})^2}  \mathcal{M}^{\al,p}(\gamma) \rho(x)\rho(y) dx dy.
\end{align}

\subsection{Known results}

Various discretization of O'Hara energy has been considered mainly due to the purpose of numerical calculation.

First, D. Kim and R. Kusner introduced the discretization of M\"{o}bius energy for polygons in \cite{KK}.

This energy, defined on the class of arc length parametrizations of polygons of length $\mathcal{L}$ with $n$ segments, is given by
\[\mathcal{E}_n(P):=\sum_{\substack{i,j=1 \\ i\neq j}}^n \left(\frac{1}{|P(a_j)-P(a_i)|^2}-\frac{1}{d(a_j,a_i)^2} \right) d(a_{i+1},a_i)d(a_{j+1},a_j),\]
where the $a_i$ are consecutive points on $\TL$, or interval $[0,\mathcal{L}]$ if we consider the polygon parametrised over an interval. This energy is scale invariant. A slight variant would be to
take $2^{-1}(d(a_{k-1},a_k)+d(a_k,a_{k+1}))$ instead of $d(a_{k+1}, a_k)$.

S. Scholtes proved that this discretization $\Gamma$-converges to the M\"{o}bius
energy in \cite{Sch}. He furthermore showed that this energy is minimized by regular n-gons.

Second, J. Simon \cite{JS} defined the so-called {\it{minimal distance energy}} for a polygon $P$ by

\[\mathcal{E}^m_s(P)=\tilde{\mathcal{E}}^m_s(P)-\tilde{\mathcal{E}}^m_s(R_m)+4\] with \[\tilde{\mathcal{E}}^m_s(P)=\sum_{|i-j|>1} \frac{|X_i||X_j|}{\mathrm{dist}(X_i,X_j)^2} ,\] 
where $R_n$ is the regular $n$-gon. Note, that this energy is scale invariant. 
Third, in \cite{Nag} is established M\"{o}bius invariant discrete energy and show $\Gamma^{-}$-convergence in $W^{1,q}$-metric sense. The definition of that energy is as follows:
\[ \mathcal{E}^m_{\cos}(P)= \sum_{d_m(i,j)>1} \frac{|\Delta_i P| |\Delta_j P|}{|\Delta^j_i P||\Delta^{j+1}_{i+1} P|} \left(1- \frac{1}{2} (\cos(\al_{i j})+ \cos(\tilde{\al}_{i j}) )\right), \]
where $P(\theta_i)$ is a vertices of closed polygon, $i = 1,2,...,m$ and  
$\Delta^{j}_i P:= P(\theta_j)-P(\theta_i)$, $\Delta_i P= \Delta^{i+1}_{i} P$,
$\al_{i j}$ be the angle of crossing of the circles through at the points $P(\theta_i),P(\theta_{i+1}),P(\theta_j)$ and $P(\theta_j),P(\theta_{j+1}),P(\theta_i)$ and $\tilde{\al}_{i j}$ be the angle of crossing of the circles through at the points $P(\theta_i),P(\theta_{i+1}),P(\theta_{j+1})$ and $P(\theta_j),P(\theta_{i+1}),P(\theta_{j+1})$.

$\newline$
Our result is to construct random discretization of O'Hara energy, and to show $\Gamma$-convergence in that space based on optimal transport theory. We even show {\bf{locally uniform convergence}} and {\bf{compactness}}.

\subsection{Main results}

We defined new discretization of O'Hara $(\alpha,p)$-energy using random variable on $\TL$.

\begin{Dfn}[Random O'Hara Energy]

Let $\left\{X_i \right\}_{i \in \N}$ be a sequence of i.i.d. random variable on $\TL$ with probability density function $\rho$.

Random O'hara energy $R_{n,\rho} \mathcal{E}^{\alpha,p}:\mathcal{A} \rightarrow \R \cup \{+\infty\}$ is defined as follow:
\[R_{n,\rho}\mathcal{E}^{\al,p} (\gamma):=\frac{1}{n^2} \sum_{\substack{i,j=1 \\ i\neq j}}^n \left( \frac{1}{|\gamma(X_i)-\gamma(X_j) |^\al}-\frac{1}{\mathcal{D}(\gamma(X_i),\gamma(X_j))^\al} \right)^p.\]
\end{Dfn}
\begin{Rem}
Since $\{X_i \}_{i\in\N} $ has probability density function $\rho$, we always have 

$\mathbb{P}(X_i = X_j)=0,$ for any  $i \neq j$. Therefore we can defined $R_{n,\rho}\mathcal{E}^{\al,p}$ almost surely.

\end{Rem}

Then we introduce the space for comparing continuous model and discrete model as follows.

\begin{Dfn}[The $TL^q$ metric space \cite{Garcia1}] $TL^q$ metric is defined on particular spaces of the family
\[TL^q(\mathbb{R} / \mathcal{L}\mathbb{N}):=\left\{(\mu,f) \mid \mu \in \mathcal{P}(\mathbb{R}/\mathcal{L}\mathbb{N}),f\in L^q(\mathbb{R}/\mathcal{L}\mathbb{Z};\mu) \right\}, \] 
where $1\leq q <\infty$ and $\mathcal{P}(\mathbb{R}/\mathcal{L}\mathbb{N})$ denotes the set of Borel probability measure on $\mathbb{R}/\mathcal{L}\mathbb{N}.$ For $(\mu,f)$ and $(\nu,g)$ in $TL^q$ we define the distance
\[d_{TL^q}\bigl((\mu,f),(\nu,g)\bigr)
:= \inf_{\pi \in \Gamma(\mu,\nu)} \left(\int_{(\mathbb{R}/\mathcal{L}\mathbb{Z})^2} \left(|x-y|^q+|f(x)-g(y)|^q \right) d\pi(x,y) \right)^{1/q},\]
where $\Gamma(\mu,\nu)$ is the set of all coupling (or transportation plans) between $\mu$ and $\nu$, that is, the set of all Borel probability measures on $(\TL)^2$ for which the marginal on the first variable is $\mu$ and the marginal on the second variable is $\nu$.

It is shown in \cite{Garcia2} the $d_{TL^q}$ is actually a metric. 

%As discussed in 
The distance $d_{TL^q}$ called a transportation distance between functions defined on graph. The $TL^q$ topology provides a general and versatile way to compare functions in a discrete setting with functions in a continuum setting. It is a generalization of the weak convergence of measures and of $L^q$ convergence of functions.

\end{Dfn}

\begin{Dfn} Let $\{X_i\}_{i\in \N}$ be a sequence of i.i.d. random variable and let us denote by $\nu_n$ the {\it{empirical measure}} of $\{X_i\}_{i\in\N}$:
\[\nu_n:=\frac{1}{n} \sum_{i=1}^n \delta_{X_i},\]
where $\delta_{X}$ is Dirac measure of $X$.
\end{Dfn}

\begin{Dfn} 
Let $\{X_i\}_{i \in \N}$ be a sequence of i.i.d. random variable on $\TL$ and $\nu_n$ be a empirical measure of $\{X_i\}_{i\in \N}$ and $\nu$ be a distribution measure of $\{X_i\}_{i \in \N}$ then, we use a slight abuse of notation and write $(\nu_n,\gamma_n) \xrightarrow{TL^q} (\nu ,\gamma)$ as $\gamma_n \xrightarrow{TL^q} \gamma.$   

\end{Dfn}

The main results of the paper is
\begin{Thm}[$\Gamma$-convergence and locally uniform convergence]
Let $\{ X_i \}_{i\in\N} $ be a sequence of i.i.d. random variable with probability density function $\rho$ on $\TL$.  
 Then $R_{n,\rho}\mathcal{E}^{\al,p}$ $\Gamma$-converge to $\mathcal{E}^{\al,p}_{\rho}$ as $n \rightarrow \infty$
in the $TL^1$ sense.

Moreover, we set $\mathcal{F}:=\{\gamma \in \mathcal{A} \mid \mathcal{E}^{\al,p}_{\rho}(\gamma) < \infty \}$, then $R_{n,\rho}\mathcal{E}^{\al,p}\lfloor_{\mathcal{F}}$ locally uniformly converge to $\mathcal{E}^{\al,p}_{\rho}\lfloor_{\mathcal{F}}$ as $n \rightarrow \infty$
in the $TL^1$ sense a.s. $\omega \in \Omega$.

\end{Thm}

\begin{Thm}[Compactness]
Let $\rho$ be a bounded from below by a positive constant and let $2 \leq \al p < 2p+1$ and $1 \leq q < \infty$,

Assume $\{\gamma_n\}_{n\in \mathbb{N}} \subset TL^q(\TL)$ satisfying 
\[\sup_{n \in \N}  \left( R_{n,\rho} \mathcal{E}^{\al,p}(\gamma_n) + \|\gamma_n \|_{L^1(\T,\nu_n)} \right)  < \infty. \]
Then $ \{\gamma_n \}_{n\in \N}$ is relatively compact in the $TL^q(\TL)$ sense a.s. $\omega \in \Omega$.

\end{Thm}

The metric used in the $\Gamma$-convergence and locally uniform convergence is the $TL^1$ sense, which will be defined in Section 2.

\section{Preliminaries}

\subsection{$\Gamma$-convergence on metric spaces and locally uniformly convergence in the $TL^q$ sense}

We recall notation of general $\Gamma$-converge and locally uniform convergence in the $TL^q$ sense. 

\begin{Dfn}[$\Gamma$-convergence on metric spaces] Let $(X,d)$ be a metric space. Let $F_n:X\rightarrow [0,\infty]$ be a sequence of functionals. The sequence $\{F_n\}_{n\in \N}$ $\Gamma$-converges with respect to metric $d$ to the functional $F:X\rightarrow [0,\infty]$ as $n\rightarrow \infty$ if the following inequality hold:

\begin{itemize}
\item[i)] For every $x\in X$ and every sequence $\{x_n\}_{n\in \N}$ converging to $x$
\[\liminf_{n\rightarrow \infty} F_n(x_n) \geq F(x),\]
\item[ii)] For every $x\in X$ there exists a sequence $\{x_n\}_{n\in N}$ converging to $x$ satisfying
\[\limsup_{n\rightarrow \infty}F_n(x_n) \leq F(x).\]
\end{itemize}

\end{Dfn}

\begin{Dfn}[Locally uniform convergence in the $TL^q$ sense]
A set $K \subset \mathcal{A}$ {\it{is sequentially compact in the $TL^q$ sense}} if it satisfies following conditions.

For all sequence $\{\gamma_n\}_{n\in \N} \subset K$ , there is a subsequence $\{\gamma_{n_k}\}_{k\in \N}$ and a $\gamma \in K$ such that $\gamma_{n_k} \xrightarrow{TL^q} \gamma$ as $k\rightarrow \infty$.

Let $X$ be a space containing $L^q(\nu_n)$ and $L^q(\nu)$,

and let $F_n:X \rightarrow \R$ and $F:X \rightarrow \R.$

The sequence $\{F_n\}_{n \in \N}$ {\it{locally uniformly converges to $F$ in the $TL^q$ sense}} if it satisfies following conditions.

For any sequentially compact set $K \subset X$ in the $TL^q$ sense,
\[\lim_{n\rightarrow \infty} \sup_{\gamma \in K} |F_n(\gamma)-F(\gamma)|=0.\]

\end{Dfn}

$\newline$
We first discuss an equivalent condition for locally uniform convergence in the $TL^q$ sense. 

\begin{Prop} 

Let $F:X\rightarrow \R$ is continuous and a sequence $\{F_n\}_{n \in \N}$ locally uniformly converge to $F$ in the $TL^q$ sense.

if and only if

For any sequence $\{\gamma_n\}_{n\in \N} \subset X $ with $\gamma_n \xrightarrow{TL^q} \gamma$,  
\[\lim_{n\rightarrow \infty} F_n(\gamma_n)=F(\gamma).\]
\end{Prop}

This can be proved in a similar way to prove Ascoli-Arzelà theorem \cite[Theorem 7.25.]{Rudin}.
\begin{proof}

For $\gamma\in \mathcal{A}$ and $r>0$, we set $\D \mathcal{B}(\gamma,r):=\{\tilde{\gamma}\in X \mid \mbox{There exists an } n \in \N$ such that $d_{TL^q}\bigl((\nu_n,\gamma),(\nu,\tilde{\gamma}) \bigr) < r \}.$

First, we show that suppose $K\subset X$ be a sequentially compact set in the $TL^q$ sense, then for any $\ep>0$, there is a sequence $\{\gamma_i\}_{i=1}^{N_{\ep}} \subset K,$ such that $\D \bigcup_{i=1}^{N_{\ep}} \mathcal{B}(\gamma_i,\ep) \supset K.$ 

If not, there is a $r>0$ such that for all $\{\gamma_i\}_{i=1}^m \subset K$, $\D K \setminus \bigcup_{i=1}^m \mathcal{B}(\gamma_i,r) \neq \emptyset.$

We choose $\gamma_1 \in K$ and inductively choose $\D \gamma_n \in K \setminus \bigcup_{i=1}^{n-1} \mathcal{B}(\gamma_i,r),$ then for all $m,n\in \N,$ 
\begin{align*}
\D r &\leq d_{TL^q} \bigl((\nu_n,\gamma_n),(\nu,\gamma_m)\bigr) \leq d_{TL^q} \bigl((\nu_n,\gamma_n),(\nu_m,\gamma_m)\bigr)+d_{TL^q} \bigl((\nu_m,\gamma_m),(\nu,\gamma_m)\bigr) \\
&=d_{TL^q} \bigl((\nu_n,\gamma_n),(\nu_m,\gamma_m)\bigr).
\end{align*}

This is a contradiction to the fact that $K$ is sequentially compact set in the $TL^q$ sense.

Let $\ep_m>0$ with $\ep_m \searrow 0$. and we set 
\begin{align}
K_0:=\left\{ \{\gamma_i^{m}\}_{m=1}^{N_m} \in K \relmiddle| \bigcup_{i=1}^{N_m} \mathcal{B}(\gamma_i^{m},\ep_m) \supset K \right\}.
\end{align}
 
$\newline$
Second, we show that for all $\ep>0$, there exists a $\delta>0$ such that for all $\gamma_1,\gamma_2 \in K$ and for all $n\in \N$, if $\D d_{TL^q}\bigl((\nu_n,\gamma_1),(\nu,\gamma_2)\bigr)< \delta$ then 
\begin{align}
|F_n(\gamma_1)-F_n(\gamma_2)| &< \ep.
\end{align}
If not, there exists an $\ep>0$ such that for all $n\in \N$, there exists $\gamma_1^n,\gamma_2^n \in K$ and $m_n \in \N$
 such that $\D d_{TL^q} \bigl((\nu_n,\gamma_1^n),(\nu,\gamma_2^n) \bigr)< 1/n$ and $\D |F_{m_n}(\gamma_1^n)-F_{m_n}(\gamma_2^n)| \geq \ep.$

By the $K$ is sequentially compact, there exists a subsequence $\{\gamma_1^{n_k}\}_{k \in \N} , \{\gamma_2^{n_k}\}_{k \in \N}$ and $\gamma \in K$ such that $\gamma_1^{n_k} \xrightarrow{TL^q} \gamma$ , $\gamma_2^{n_k} \xrightarrow{TL^q} \gamma.$
Then\begin{align*}
\ep &<|F_{m_{n_k}}(\gamma_1^{n_k})-F_{m_{n_k}}(\gamma_2^{n_k})| \leq |F_{m_{n_k}}(\gamma_1^{n_k})-F(\gamma)|+|F(\gamma)-F_{m_{n_k}}(\gamma_2^{n_k})| \xrightarrow{k \rightarrow \infty} 0.
\end{align*}

This is a contradiction.

$\newline$
Third, we show that there exists a subsequence $\{F_{n_k}\}_{k \in \N}$ such that for all $j$ , $\{F_{n_k}(\gamma_j)\}_{k\in \N}$ is convergence sequence by a diagonal argument.

By $\D \lim_{n\rightarrow \infty}F_n(\gamma_1)=F(\gamma_1)$, $\{F_n(\gamma_1)\}_{n\in N} $ is bounded sequence on $\R.$

Therefore there exists a subsequence $\{F_{n(1,k)}\}_{k\in \N}$ such that $\{F_{n(1,k)}(\gamma_1)\}_{k\in \N}$ is convergence sequence on $\R.$

In the same way, there exists a subsequence $\{F_{n(2,k)}\}_{k\in \N}$ such that $\{F_{n(2,k)}(\gamma_2)\}_{k\in \N}$ is convergence sequence on $\R.$

Further in the same way, we construct subsequence $\{F_{n(p,k)}\}_{k\in \N} , p=3,4,...,$ and we set $F_{n_k}=F_{n(k,k)}$.

Finally, let any $\eta>0$, for sufficient large $m$ such that for all $n\in \N$, if $\D d_{TL^q}\bigl((\nu_n,\gamma_1),(\nu,\gamma_2)\bigr) < \ep_m$ then 
\begin{align}
|F_n(\gamma)-F_n(\gamma_i^{m})| < \eta / 3.
\end{align}
Since $\{F_{n_k}(\gamma_i^{m})\}_{k \in \N}$ is a convergence sequence on $\R$, there exists a number $N$ such that if $k,l > N$ then $|F_{n_k}(\gamma_i^{m})-F_{n_l}(\gamma_i^{m})| < \eta/2$ for $i=1,2,...,N_m$. 

Now, let $\gamma \in K$, by (5) there exist an $i$ and $n$ such that 
\[d_{TL^q}\bigl((\nu_n,\gamma_i^{m}),(\nu,\gamma)\bigr) < \ep_m.\]
By (6) and (7) if $k,l>N$ then
\begin{align*}
|F_{n_k}(\gamma)-F_{n_l}(\gamma)| &\leq |F_{n_k}(\gamma)-F_{n_k}(\gamma_i^{m})|+|F_{n_k}(\gamma_i^{m})-F_{n_l}(\gamma_i^{m})|+|F_{n_l}(\gamma)-F_{n_l}(\gamma_i^{m})| \\
&\leq \eta.
\end{align*}
This indicates that the uniform convergence on $K$.

$\newline$
Assume that $F$ is continuous and $F_n$ locally uniformly converge to $F$ in the $TL^q$ sense and $\gamma_n \xrightarrow{TL^q} \gamma.$

Clearly $\{\gamma_n \mid n\in \N \} \cup \{\gamma \} \subset X$ is sequentially compact in the $TL^q$ sense.

Therefore
\begin{align*}
|F_n(\gamma_n)-F(\gamma)| &\leq |F_n(\gamma_n)-F(\gamma_n)|+|F(\gamma_n)-F(\gamma)| \\
&\leq \sup_{y \in\{\gamma_n\}_{n\in \N} \cup \{\gamma\}} \bigl(|F_n(y)-F(y)|\bigr)+|F(\gamma_n)-F(\gamma)| \xrightarrow{n\rightarrow \infty} 0.
\end{align*}

\end{proof}

\subsection{The property of $TL^q$ space and empirical measure}

In this subsection, we consider the space based on optimal transport theory to compare discrete and continuous model.

Given a Borel map $T:\TL \rightarrow \TL$ and $\mu \in \mathcal{P}(\TL)$ the {\it{push-forward}} of $\mu$ by $T$, denoted by $T_{\#} \mu \in \mathcal{P}(\TL)$ is given by:
\[T_{\#}\mu(A):=\mu(T^{-1}(A)) , A \in \mathcal{B}(\TL).\]

\begin{Dfn}[\cite{Garcia1}]

We say that a sequence of transportation plans $\left\{\pi_n \right\}_{n\in \mathbb{N}} \subset \Gamma(\mu,\mu_n) $ is {\it{stagnating}} if it satisfies 
\[\lim_{n \rightarrow \infty} \int_{(\TL)^2} |x-y|^q d\pi_n (x,y)=0.\]
$\mu,\mu_n \in \mathcal{P}(\TL),T_n:\TL\rightarrow \TL :$ transportation maps with ${T_n}_{\#} \mu=\mu_n $.

We say that a sequence of transportation maps $\{T_n\}_{n\in \mathbb{N}}$ is {\it{stagnating}} if it satisfies
\[\lim_{n \rightarrow \infty} \int_{\TL} |x-T_n(x)|^q d\mu(x)=0.\]
\end{Dfn}

Accept the following proposition introduced by \cite{Garcia1}
\begin{Prop}[\cite{Garcia1}]
Let $(\mu,\gamma)\in TL^q$ and let $\left\{(\mu_n,\gamma_n)\right\}_{n\in \mathbb{N}} \subset TL^q $．

The following statements are equivalent;
\begin{itemize}
\item[i)] $(\mu_n,\gamma_n) \rightarrow (\mu,\gamma)$ in the $TL^q$.

\item[ii)] $\mu_n \rightharpoonup \mu $ and  for every stagnating sequence of transportation plans $\left\{\pi_n \right\}_{n\in \N} \subset \Gamma(\mu,\mu_n)$ 
\begin{align}
\int_{(\TL)^2} |\gamma(x)-\gamma_n(x)|^q d \pi_n(x,y) \rightarrow 0.
\end{align}
\item[iii)] $\mu_n \rightharpoonup \mu $ and there exists a stagnating sequence of transportation plans  $\left\{\pi_n \right\}_{n\in \mathbb{N}} \subset \Gamma(\mu,\mu_n)$ such that
\begin{align}
\D \int_{(\TL)^2} |\gamma(x)-\gamma_n(x)|^q d \pi_n(x,y) \rightarrow 0.
\end{align}
Moreover, if the measure $\mu$ is absolutely continuous with respect to the Lebesgue measure, the following are equivalent to the previous statement:

\item[iv)] $\mu_n \rightharpoonup \mu$ and there exists a stagnating sequence of transportation maps $\left\{ T_n \right\}_{n \in \N}$ $($ with ${T_n}_{\#}\mu=\mu_n  )$ such that
\begin{align}
\D \int_{\TL} |\gamma(x)-\gamma_n(T_n(x))|^q d\mu(x) \rightarrow 0. 
\end{align}
\item[v)] $\mu_n \rightharpoonup \mu$ and for any stagnating sequence of transportation maps $\left\{ T_n \right\}_{n \in \N}$ with ${T_n}_{\#}\mu=\mu_n  )$， 
\begin{align}
\int_{\TL} |\gamma(x)-\gamma_n(T_n(x))|^q d\mu(x) \rightarrow 0.
\end{align}

\end{itemize}
\end{Prop}

\begin{Rem} Thanks to Proposition 2.2 when $\mu$ is absolutely continuous with respect to the Lebesgue measure $\gamma_n \xrightarrow{TL^p} \gamma$ as $n \rightarrow \infty$ if and only if for every (or one) $\{T_n \}_{n\in \N}$ stagnating sequence of transportation maps (with ${T_n}_{\#} \mu=\mu_n $)  $\gamma_n \circ T_n \xrightarrow{L^p(\mu)} \gamma $ as $n \rightarrow \infty.$ %(this in particular implies the last part of あ). 
 Also $\{u_n\}_{n \in \N} $ is relatively compact in $TL^p$ if and only if for every (or one) $\{T_n \}_{n \in N}$ stagnating sequence of transportation maps (with $T_{\#} \mu=\mu_n $) $\{u_n \circ T_n \}_{n \in \N}$ is relatively compact in $L^p(\mu)$.

\end{Rem}

We recall the following proposition.
\begin{Prop}[{\bf{Glivenko-Cantelli's Theorem}}]
Let $\{X_i\}_{i \in \N} $ be a sequence of i.i.d. random variable on $\TL$, and let $\nu$ is distribution measure of $\{X_i\}_{i \in \N}$, and $\nu_n$ is empirical measure of $\{X_i\}_{i \in \N}$.

$F_n(x):=\nu_n((-\infty,x])$ and $F$ is distribution function of $\{X_i\}_{i \in \N}$, then,
\[\lim_{n \rightarrow \infty} \|F_n-F\|_{\infty}=0 \qquad a.s. ~ \omega \in \Omega. \]

\end{Prop}

\begin{Thm}[\cite{Garcia3}] Let $D \subset \R^d$ be a bounded, connected, open set with Lipschitz boundary. 
Let $\nu$ be a probability measure on $D$ with density $\rho:D\rightarrow (0,\infty)$ {which is bounded from below and from above by positive constants.} Let $\{X_i\}_{i\in \N}$ be a sequence of independent random {points} distributed on $D$ according to measure $\nu$ and let $\nu_n$ {be the associated empirical measures.}  Then there is a constant $C>0$ such that for $a.s.~ \omega \in \Omega$ there exists a sequence of transportation maps $\{T_n\}_{n\in \N}$ from $\nu$ to $\nu_n$ $({T_n}_{\#} \nu=\nu_n)$ and such that

if $d=2$ then \[\limsup_{n\rightarrow \infty} \frac{n^{1/2} \|Id-T_n\|_{\infty}}{(\log n)^{3/4}} \leq C\]

and if $d \geq 3$ then \[\limsup_{n\rightarrow \infty} \frac{n^{1/d} \|Id-T_n\|_{\infty}}{(\log n)^{1/d}} \leq C.\]

\end{Thm}

\section{$\Gamma$-convergence and locally uniformly convergence}

{\bf{Proof of Theorem 1.1}}

\begin{proof}

Let $\nu_n$ be an empirical measure of $\mathcal{L}^1 \lfloor \rho$ and $T_n:\TL \rightarrow \TL$ be a transportation maps with ${T_n}_{\#}\mathcal{L}^1 \lfloor \rho=\nu_n$.

・{\bf{liminf inequality}}

Assume that $\gamma_n \rightarrow \gamma$ in $TL^1 $ as $n \rightarrow \infty$. Since ${T_n}_{\#}\mathcal{L}^1 \lfloor \rho=\nu_n$, we change variables to get

\begin{align}
\D R_{n,\rho}\mathcal{E}^{\al,p}(\gamma_n)&=\int_{(\TL)^2} \mathcal{M}^{\alpha,p}(\gamma_n) d\nu_n(x) d\nu_n(y) \\
&=\int_{(\TL)^2} \mathcal{M}^{\alpha,p}(\gamma_n \circ T_n) \rho(x) \rho(y) dx dy.
\end{align}
By the way we notice that
\[|\gamma_n(T_n(x))-\gamma_n(T_n(y)) | \leq |\gamma_n(T_n(x))-\gamma(x)|+|\gamma(x)-\gamma(y)|+|\gamma(y)-\gamma_n(T_n(y)) | .\]
and 
\begin{align*}
\D \mathcal{D}(\gamma(x),\gamma(y)) &\leq \mathcal{D}(\gamma(x),\gamma_n(T_n(x)))+\mathcal{D}(\gamma_n(T_n(x)),\gamma_n(T_n(y)))+\mathcal{D}(\gamma_n(T_n(y)),\gamma(y)) \\
&\leq |x-T_n(x)|+\mathcal{D}(\gamma_n(T_n(x)),\gamma_n(T_n(y)))+|T_n(y)-y| \\
&\leq 2 \|T_n-\mathrm{Id}\|_{\infty}+\mathcal{D}(\gamma_n(T_n(x)),\gamma_n(T_n(y))).
\end{align*}

By Proposition 2.2 we deduce $\gamma_n \circ T_n \rightarrow \gamma $ in $L^1(\TL) $.

Thus, by taking an appropriate subsequence $\left\{\gamma_{n_k} \circ T_{n_k} \right\}_{k\in \mathbb{N}} $ of $\left\{\gamma_n \circ T_n \right\}_{n\in \mathbb{N}}$, we deduce
$\gamma_{n_k}(T_{n_k}(x)) \rightarrow \gamma(x) \ \  \mbox{a.e.}~x \in \TL$
and therefore
\[ \liminf_{n\rightarrow \infty} \mathcal{M}^{\al,p}(\gamma_n \circ T_n)  \geq  \mathcal{M}^{\al,p}(\gamma) \ \ \mbox{a.e.}~(x,y)\in (\TL)^2.\]
So that $\displaystyle \mathcal{M}^{\al,p}(\gamma_n \circ T_n) > 0$, using Fatou's lemma we get
\[ \liminf_{n \rightarrow \infty} R_{n,\rho}\mathcal{E}^{\al,p}_n(\gamma_n) \geq \mathcal{E}^{\al,p}_{\rho}(\gamma).\]
・{\bf{limsup inequality}}

Assume that $\gamma_n\rightarrow \gamma \ \ \mathrm{in} \ \ TL^1 $ as $n\rightarrow \infty$.

If $\D \mathcal{E}^{\alpha,p}_{\rho}(\gamma)=\infty ,$ we are done, and so we henceforth assume $\D \mathcal{E}^{\alpha,p}_{\rho}(\gamma) < \infty.$ 

We observe that
\begin{align*}
|\gamma(x)-\gamma(y) | \leq |\gamma(x)-\gamma_n(T_n(x)) |+ |\gamma_n(T_n(x))-\gamma_n(T_n(y)) |+ |\gamma_n(T_n(y))-\gamma(y) |,
\end{align*}
and
\begin{align*}
\mathcal{D}(\gamma_n(T_n(x)),\gamma_n(T_n(y))) \leq \mathcal{D}(\gamma_n(T_n(x)),\gamma(x))+\mathcal{D}(\gamma(x),\gamma(y))+\mathcal{D}(\gamma(y),\gamma_n(T_n(y))). 
\end{align*}

In the same way as liminf  inequality, we get 

$\D \limsup_{n \rightarrow \infty} \mathcal{M}^{\al,p}(\gamma_n \circ T_n)  \leq \mathcal{M}^{\al,p}(\gamma) \ \ \mbox{a.e.}~(x,y)\in (\mathbb{R}/\mathbb{Z})^2.$

Since $\mathcal{E}^{\al,p}_{\rho}(\gamma) < \infty$, using Fatou's lemma to $ \mathcal{M}^{\al,p}(\gamma)-\mathcal{M}^{\al,p}(\gamma_n \circ T_n) ,$ we get
\[ \limsup_{n \rightarrow \infty} R_{n,\rho}\mathcal{E}^{\al,p}(\gamma_n) \leq \mathcal{E}^{\al,p}_{\rho}(\gamma).\]
By Proposition 2.1, the proof is now complete.

\end{proof}

\section{Compactness}

In this section, we would like to prove Theorem 1.2 we first recall several function space.
 
\subsection{Function spaces} 
\begin{Dfn}[{\bf{Sobolev-Slobodeckij spaces}}]

For $s\in (0,1)$ and $q\in [1,\infty)$ we set
\[[\gamma]_{W^{s,p}} := \left( \int_{\TL} \int_{-\mathcal{L}/2}^{\mathcal{L}/2}  \frac{|\gamma(u+w)-\gamma(u)|^p}{|w|^{1+p s}} dw du \right)^{1/p}\]
\[W^{s,p}(\TL,\mathbb{R}^n):= \left\{\gamma \in L^q(\TL,\mathbb{R}^n) \mid [\gamma]_{W^{s,p}} < \infty \right\}\]
and equip this space with the norm
\[\|\gamma \|_{W^{s,q}(\mathbb{R}/\mathbb{Z},\mathbb{R}^n)}:=\|\gamma \|_{L^q}+[\gamma]_{W^{s,p}}.\]
Furthermore, we let
\[W^{1+s,q}(\mathbb{R}/\mathbb{Z},\mathbb{R}^n):= \left\{ \gamma \in W^{1,q}(\mathbb{R}/\mathbb{Z},\mathbb{R}^n) \mid \gamma ' \in W^{s,q}(\mathbb{R}/\mathbb{Z},\mathbb{R}^n) \right\}\]
and\[\|\gamma \|_{W^{1+s,q}}:=\left(\|\gamma \|^{q}_{W^{s,q}}+\|\gamma \|^{q}_{W^{1,q}} \right)^{1/q}.\]
\end{Dfn}

\begin{Dfn}[{\bf{Besov spaces}}] For $s\in \R $, $1 \leq p,q \leq \infty$ and $\mathcal{S}$ is the Schwartz space.

We set\[\psi(\xi)= \begin{cases}
    1 & (|\xi| \leq 4) \\
    0 & (|\xi| \geq 8)
  \end{cases} , ~~\varphi(\xi)= \begin{cases}
    1 & (2 \leq |\xi| \leq 4) \\
    0 & (|\xi| \leq 1 \mbox{ or } |\xi| \geq 8 )
  \end{cases},\]
  and we set $\varphi_j(\xi)=\varphi(2^{-j} \xi)$ and $\tau(D)f:=\mathcal{F}^{-1}[\tau \cdot \mathcal{F}f]$ for $\tau \in \mathcal{S}$ and $f\in \mathcal{S}'$ we set
  
  \[\|f\|_{B^s_{p,q}}:= \begin{cases}
    \D \|\psi(D) f\|_{L^p}+ \left( \sum_{j=1}^{\infty} 2^{jqs} \|\varphi_j(D)f \|_{L^p}^q \right)^{1/q} & (1\leq q < \infty) \\
    \D \|\psi(D) f\|_{L^p}+\sup_{j\in \N} 2^{js} \|\varphi_j(D)f\|_{L^p} & (q=\infty)
  \end{cases} \]
  \[B^s_{p,q}:=\{f\in \mathcal{S}' \mid \|f\|_{B^s_{p,q}} < \infty \}.\]
\end{Dfn}

Note that $W^{s,p}$ agrees with $B^s_{p,p}$.

We recall the following theorem.

\begin{Prop}[{\bf{Embedding Besov spaces}}] 
Let $0<p_0,p_1\leq \infty$ and $0< q_0,q_1 \leq \infty$ and $-\infty< s_1 < s_0 < \infty$.
Assume that $s_0-\frac{n}{p_0} > s_1-\frac{n}{p_1}$, then
\[B^{s_0}_{p_0,q_0} \hookrightarrow B^{s_1}_{p_1,q_1}.\]

\end{Prop}

\begin{Thm}[{\bf{Kondrachov embedding theorem}}]
Let $1 \leq p,q < \infty$ and $1 < k,s $.

Assume that $k-\frac{1}{p}>s-\frac{1}{q}$, then the Sobolev embedding 
\[W^{k,p}(\TL) \hookrightarrow  W^{s,q}(\TL)\]
is completely continuous.

\end{Thm}

\begin{Thm}[{\bf{Gagliardo-Nirenberg interpolation inequality}}]

Suppose that $f \in W^{l,q}(\TL) \cap W^{j,r}(\TL) $ for some $j<l , q,r \geq 1.$

Assume that $j\leq k < l$ and for some $0<\Theta \leq 1$,

$k-\frac{n}{p}=\Theta \left(l-\frac{n}{q} \right)+(1-\Theta) \left( j-\frac{n}{r} \right),$ then $f \in W^{k,p}(\TL)$. Moreover, we have for all $f \in W^{l,q}(\TL) \cap W^{j,r}(\TL)$,

\[\|f\|_{W^{k,p}}\leq C \|f\|^{\Theta}_{W^{l,q}} \|f\|^{1-\Theta}_{W^{j,r}} \]

for some constant $C>0$ independent of $f$.
\end{Thm}

\begin{Prop}[{\bf{Embedding $L^1$ space for Besov space}}]
Let $s\in \R$ and $0 < q \leq \infty$, then
\[L^1 \hookrightarrow B^{0}_{1,q} \mbox{ if and only if } q= \infty.\]

\end{Prop}

\begin{Thm} Let $-\infty < s_0 < \infty$ , $-\infty < s_1 < \infty$ , $0 < p_0 \leq \infty$ , $0 < p_1 \leq \infty$ , $0 < q_0 \leq \infty$ , $s= (1-\Theta) s_0+\Theta s_1$. Assume that for some $1< \Theta \leq 1$, $ \frac{1}{p}=\frac{1-\Theta}{p_0}+\frac{\Theta}{p_1}$ and $ \frac{1}{q}=\frac{1-\Theta}{q_0}+\frac{\Theta}{q_1},$ then
\[(B^{s_0}_{p_0,q_0}(\TL),B^{s_1}_{p_1,q_1}(\TL))_{\Theta,p}=B^{s}_{p,q}(\TL).\]

\end{Thm}

The following theorem is based on \cite{Blatt}, and explain conditions for which O'Hara energy becomes finite.

\begin{Thm}[\cite{Blatt}] Let $\gamma \in \mathcal{A}$ and $\al,p \in (0, \infty)$ with $\al p \geq 2$ and $s:=\frac{\al p-1}{2p} <1 $ and $p \geq 1$, then $\mathcal{E}^{\al,p}(\gamma) < \infty $ if and only if $\gamma \in W^{1+s,2p}(\TL,\mathbb{R}^n).$ Moreover, there is a $C = C(\al,p)$ such that 
\[\|\gamma ' \|_{W^{s,2p}}^{2p} \leq C(\mathcal{E}^{\al,p}(\gamma)+\|\gamma ' \|_{L^{2p}}^{2p}).\]

\end{Thm}

\subsection{Proof of Theorem 1.2.}

\begin{proof}
Let $ s:=\frac{\al p -1}{2p}$. By Theorem 4.4 we see
\[\|\gamma \|^{2p}_{W^{1+s,2p}} \leq C \left(\mathcal{E}^{\al,p}(\gamma)+\|\gamma \|^{2p}_{W^{1,2p}} \right) .\]
By Theorem 4.3 we choose $\Theta >0$ such that $ \Theta < \frac{(2p-1)q+2p}{(2p +\al p-2)q+2p} \leq 1$, and $t$ with $ \frac{1}{2p}-1 > t$， to get
\begin{align}
\D \|\gamma \|_{W^{1,2p}} &\leq \|\gamma\|_{W^{1+s/2,2p}} \\
&=\|\gamma \|_{B^{1+s/2}_{2p,2p}} \leq C \|\gamma\|^{\Theta}_{B^{1+s}_{2p,2p}} \|\gamma \|^{1-\Theta}_{B^t_{2p,2p}}.
\end{align}
Since $L^1 \hookrightarrow B^{0}_{1,\infty} \hookrightarrow B^{t}_{2p,2p}$ this implies yields 
\[\D \|\gamma \|_{W^{1,2p}}  \leq C \|\gamma\|^{\Theta}_{B^{1+s}_{2p,2p}} \|\gamma \|^{1-\Theta}_{L^1}.\]
Using Young's inequality, for all $ \varepsilon >0$,  we get
\begin{align}
\D \|\gamma \|^{2p}_{W^{1,2p}} &\leq C \|\gamma \|^{2p \theta}_{W^{1+s,2p}} \|\gamma \|^{2p(1-\Theta)}_{L^1} \\
&\leq C \ep^{1/\Theta}  \Theta \|\gamma \|^{2p}_{W^{1+s,2p}} + C(1-\Theta)\frac{\|\gamma \|^{2p}_{L^{1}}}{\ep^{1/(1-\Theta)} }.
\end{align}
Therefore, for sufficient small $\ep>0$, we conclude that
\begin{align}
\|\gamma \|^{2p}_{W^{1+s,2p}} \leq C' \left( \mathcal{E}^{\al,p}(\gamma)+ \|\gamma \|^{2p}_{L^1} \right).
\end{align}

Let $\{\gamma_n\}_{n \in \N}$ be a sequence of $TL^q(\T)$ with 
\[\D \sup_{n \in \N} \left( R_{n,\rho} \mathcal{E}^{\al,p}(\gamma_n)+\|\gamma_n\|_{L^1(\TL,\nu_n)} \right)< \infty ,\]
and let $T_n:\TL \rightarrow \TL$ be a transportation maps with ${T_n}_{\#} \mu = \mu_n $, then
\[\sup_{n \in \N} \mathcal{E}^{\al,p}_{\rho} (\gamma_n \circ T_n) < \infty.\] 
Since $\rho$ is bounded from below by a positive constant, we deduce that 
\[\sup_{n \in \N} \mathcal{E}^{\al,p} (\gamma_n \circ T_n) < \infty.\]
Therefore by (18), we see 
\[\sup_{n\in \N} \|\gamma_n \circ T_n \|_{W^{1+s,2p}} < \infty.\] 
Since $ 1+s-\frac{1}{2p}=\frac{2p+\al p-2}{2p} \geq 0 > -\frac{1}{q}$, Theorem 4.1, yield a compact embedding 
\[\iota: W^{1+s,2p}(\TL) \hookrightarrow \hookrightarrow L^q(\TL).\]
Therefore there exists a $\{n_k \}_{k \in \N} \subset \N$ and $\gamma \in L^q(\TL,\rho)$
such that $\gamma_{n_k} \circ T_{n_k} \rightarrow \gamma $ in $L^q(\T,\rho).$

By Proposition 2.2 we see $\gamma_{n_k} \rightarrow \gamma $ in $TL^q(\TL)$.

\end{proof}

\end{document}